\documentclass[12pt,twoside]{amsart}
\usepackage{amssymb}
\usepackage{amscd}

\title{Three-dimensional terminal toric flips}
\author{Osamu Fujino, Hiroshi Sato, Yukishige Takano, Hokuto Uehara} 
\date{2008/4/7}
\address{Graduate School of Mathematics\\ 
 Nagoya University, Chikusa-ku Nagoya 464-8602 Japan}
\email{fujino@math.nagoya-u.ac.jp}

\address{Faculty of Economics and Information\\ 
Gifu Shotoku Gakuen University, 
1-38 Nakauzura Gifu 500-8288 Japan}
\email{hirosato@gifu.shotoku.ac.jp}

\address{Department of Mathematics\\ Tokyo Metropolitan University, 
1-1 Minami-Ohsawa Hachioji-shi Tokyo 192-0397 Japan}

\address{Department of Mathematics\\ Tokyo Metropolitan University, 
1-1 Minami-Ohsawa Hachioji-shi Tokyo 192-0397 Japan}
\email{hokuto@tmu.ac.jp}

\newcommand{\Spec}[0]{{\operatorname{Spec}}}

\newtheorem{thm}{Theorem}[section]

\newtheorem{cor}[thm]{Corollary}

\newtheorem{claim}{Claim}

\theoremstyle{definition}
\newtheorem{ca}{Case}

\newtheorem{rem}[thm]{Remark}
\newtheorem*{ack}{Acknowledgments}       
 
\newtheorem*{notation}{Notation}         

\begin{document}
\bibliographystyle{amsalpha+}

\begin{abstract}
We describe three-dimensional terminal toric flips. 
We obtain the complete local 
description of three-dimensional terminal toric flips.
\end{abstract}

\maketitle

\section{Introduction} 
The main purpose of this paper is to describe 
three-dimensional 
terminal toric flips. 

\begin{thm}[{cf.~Theorems \ref{main1} and \ref{main-taka}}]\label{11}
Let $\varphi:X\to Y$ be a small proper toric 
morphism such that $X$ is a three-dimensional 
toric variety with only terminal singularities. 
Note that $X$ is not assumed to be $\mathbb Q$-factorial. 
Let $C\simeq \mathbb P^1$ be an exceptional 
curve of $\varphi$. 
Assume that $-K_X\cdot C>0$. 
Then one of the torus invariant points of $X$ on $C$ 
is non-singular and another one 
is a terminal quotient singularity. 
In particular, $X$ is $\mathbb Q$-factorial and 
is not Gorenstein around $C$. 
\end{thm}

By this result, we obtain the complete local 
description of three-dimensional terminal toric flips (see Theorem 
\ref{main1}).  
The first author apologies for the mistake in \cite[Example 4.4.2]{spe}, 
where he claims that there exist three-dimensional 
non-$\mathbb Q$-factorial terminal toric flips. 
However, Theorem \ref{11} implies that there are no such flips. 
This paper is based on the first author's private notes and 
the third author's master's thesis \cite{takano}. 

We summarize the contents of this paper. 
In Section \ref{sec2}, we describe three-dimensional 
terminal toric singularities. 
The results are well known to 
the experts. 
Section \ref{sec3} gives the complete classification 
of three-dimensional $\mathbb Q$-factorial terminal 
toric flips. It is a supplement to \cite{kmm} 
and \cite[Example-Claim 14-2-5]{ma}. In Section \ref{sec4}, we 
prove that there are no 
three-dimensional non-$\mathbb Q$-factorial 
terminal toric flips. 
Theorem \ref{main-taka} is the main theorem 
of this paper. 
The proof depends on the results in Sections \ref{sec2} and 
\ref{sec3}. 

\begin{ack}
We would like to thank Professor Kenji Matsuki, who informed us that 
Professor Yujiro Kawamata pointed out an 
error in \cite[Remark 14-2-7]{ma}. 
The first author was partially supported 
by the Grant-in-Aid for Young Scientists (A) $\sharp$17684001 from 
JSPS. 
He was also supported by the Inamori Foundation. 
The fourth author was partially supported 
by the Grant-in-Aid for Young Scientists (B) $\sharp$17740012 from JSPS.  
\end{ack}

\begin{notation}
We will work over $\mathbb C$, the complex number field, 
throughout this paper. 
Let $v_i\in N\simeq \mathbb Z^3$ for $1\leq i\leq k$. 
Then the symbol $\langle v_{1}, v_{2}, 
\cdots, v_{k}\rangle$ denotes the cone 
$\mathbb R_{\geq 0}v_{1}+\mathbb R_{\geq 0}v_2+\cdots 
+\mathbb R_{\geq 0}v_{k}$ in $N_{\mathbb R}\simeq
N\otimes_{\mathbb Z}\mathbb R$. 
\end{notation}

\section{Three-dimensional terminal toric singularities}\label{sec2} 

In this section, we characterize non-$\mathbb Q$-factorial 
affine toric threefolds with terminal singularities. 
We will use the same notation 
as in \cite{ypg}, which is an excellent exposition on terminal 
singularities. 

Let $X$ be an affine toric threefold. 
First, let us recall the following well-known theorem of 
G.~K.~White, D.~Morrison, G.~Stevens, 
V.~Danilov, and M.~Frumkin 
(see \cite[(5.2) Theorem]{ypg}). 

\begin{thm}\label{tl}
Assume that $X$ is $\mathbb Q$-factorial. 
Then $X$ is terminal if and only if $($up to 
permutations of $(x,y,z)$ and symmetries of 
${\mu}_r$$)$ $X\simeq \mathbb C^3/ {\mu _r}$ of 
type $\frac{1}{r}(a,-a,1)$ with $a$ coprime to $r$, 
where $\mu_r$ is the cyclic group of order $r$. 
In particular, if $X$ is Gorenstein and 
terminal, then $X$ is non-singular. 
\end{thm}

Here, we prove the following well-known 
result for the reader's convenience 
(cf.~\cite{ishida}, and 
\cite[Theorem 3.6]{ishi-iwa}). 

\begin{thm}\label{main}
Assume that $X$ is not $\mathbb Q$-factorial. 
Then $X$ is terminal if and only if $X\simeq \Spec \, 
\mathbb C[x,y,z,w]/(xy-zw)$. 
We call this singularity an {\em{ordinary double point}}. 
\end{thm}

By the above theorems, we obtain the 
complete list of three-dimensional terminal toric singularities. 

\begin{rem}
Mori classified three-dimensional terminal singularities. 
For the details, see \cite[(6.1) Theorem]{ypg}. 
We do not use his classification table in this paper. 
\end{rem}

\begin{proof}[{Proof of {\em{Theorem \ref{main}}}}]
Let $N=\mathbb Z^3$ and $\Delta=\langle e_1, 
\cdots, e_k\rangle$ the cone in $N$ such that 
$X=X(\Delta)$, where each $e_i$ is primitive. 
First, we prove 

\begin{claim}\label{4}
If $X$ is non-$\mathbb Q$-factorial and terminal, 
then $k=4$. 
\end{claim}
\begin{proof}[Proof of the claim] 
It is obvious that $k\geq 4$. 
Since $X$ is $\mathbb Q$-Gorenstein, there 
is a hyperplane $H\subset N$ that contains 
every $e_i$. 
On $H\simeq \mathbb Z^2$, $e_i$s span 
two dimensional convex polygon $P$. 
By renumbering $e_i$s, we can assume that they 
are arranged counter-clockwise. 
Since $X(\Delta)$ is terminal, all the lattice points 
in $P$ are $e_i$s. 
In particular, the triangle on $H$ spanned 
by $e_1$, $e_2$, and $e_3$ contains only 
three lattice points $e_i$ ($1\leq i\leq 3$) of 
$H$. 
So, after changing the coordinate of $H$, 
we can assume that $e_1=(0,1), e_2=(0,0)$, 
and $e_3=(1,0)$ in $H\simeq \mathbb Z^2$. 
It can be checked easily that 
$(1,1)\in P$ since $k\geq 4$. 
Thus, we obtain that $k=4$ and $e_4=(1,1)$. 
\end{proof}

\begin{claim}\label{cl}
Assume that $X$ is non-$\mathbb Q$-factorial, 
Gorenstein, and terminal. 
Then $X$ is isomorphic to 
$\Spec \, \mathbb C[x,y,z,w]/(xy-zw)$. 
\end{claim}
\begin{proof}[Proof of the claim]
On this assumption, the cones $\langle e_1, e_2, e_3\rangle$, 
$\langle e_1, e_2, e_4\rangle$, 
$\langle e_1, e_3, e_4\rangle$, 
and $\langle e_2, e_3, e_4\rangle$ define $\mathbb Q$-factorial 
Gorenstein affine toric threefolds with terminal 
singularities. 
By Theorem \ref{tl}, every cone listed above is non-singular. 
So, by changing the coordinate of $N$, we can assume that 
$e_1=(1,0,0)$, $e_2=(0,1,0)$, 
and $e_3=(0,0,1)$. 
Since $X$ is Gorenstein and $\langle e_1, e_2, e_4\rangle$, 
$\langle e_1, e_3, e_4\rangle$, and 
$\langle e_2, e_3, e_4\rangle$ are non-singular, 
$e_4=(-1, 1, 1)$, 
$(1, -1, 1)$, or 
$(1, 1, -1)$. 
Anyway, we can check that 
$X\simeq \Spec\, \mathbb C[x,y,z,w]/(xy-zw)$. 
\end{proof}
By the above claim, 
it is sufficient to prove 

\begin{claim}\label{6}
All the non-$\mathbb Q$-factorial toric 
affine threefolds with terminal singularities are 
Gorenstein. 
\end{claim}
\begin{proof}[Proof of the claim]  
We assume that $X$ is not Gorenstein and 
obtain a contradiction. 

Let $\overline N$ be the sublattice of $N$ spanned by 
all the lattice points on $H$ and the origin of $N$. 
In $\overline N$, 
$\Delta=\langle e_1, e_2, e_3, e_4\rangle$ defines a 
Gorenstein terminal threefold. 
So, we can assume that 
$e_1=(1,0,0)$, $e_2=(0,1,0)$, 
$e_3=(0,0,1)$, and $e_4=(1,1,-1)\in \mathbb Z^3\simeq 
\overline N$ by the proof of Claim \ref{cl}. 
First, 
we consider $\langle e_1, e_2, e_3\rangle$ in $\overline N$ 
and $N$. 
By Theorem \ref{tl}, we obtain that 
$N=\overline N +\mathbb Z\cdot \frac{1}{r}(\alpha, \beta, 
\gamma)$, where $(\alpha, \beta, 
\gamma)$ is one of the followings: 
$(a, -a, 1)$, $(a, 1, -a)$, 
$(-a, a, 1)$, $(-a,1, a)$, 
$(1, a, -a)$, $(1, -a, a)$ such 
that $0<a<r$ with $a$ coprime to $r$. 
Next, we use the terminality of 
$\langle e_1, e_2, e_4\rangle$. 
We consider the linear transform 
$T:N\to N$ such that 
$Te_1=e_1$, $Te_2=e_2$, $Te_4=e_3$.  
Then $TN=T\overline N +\mathbb Z
\cdot \frac{1}{r}(\alpha', \beta', 
\gamma')$, where $(\alpha', \beta', 
\gamma')$ is one of the followings: 
$(1+a,1-a,-1)$, $(0, 1-a, a)$, 
$(1-a, 1+a,-1)$, $(0, 1+a, -a)$, 
$(1-a, 0, a)$, $(1+a, 0, -a)$. 
Note that 
$$
\begin{pmatrix}
\alpha'\\ \beta'\\ \gamma'
\end{pmatrix}
=
\begin{pmatrix}
1 & 0 & 1 \\ 
0 & 1 & 1 \\ 
0 & 0 & -1
\end{pmatrix}
\begin{pmatrix}
\alpha\\ \beta\\ \gamma
\end{pmatrix}.
$$
We treat the first case, that is, 
$(\alpha', \beta', 
\gamma')=(1+a,1-a,-1)$. 
By the terminal lemma (see \cite[(5.4) Theorem]{ypg}), 
$r$ divides $(1+a)+(1-a)=2$ since it does not 
divide $(1+a)+(-1)$ nor $(1-a)+(-1)$. 
So, $r=2$ and $a=1$. 
Thus $\frac{1}{r}(\alpha', \beta', 
\gamma')=\frac{1}{2}(2, 0, -1)\equiv \frac{1}{2}(0,0,1)$ 
(mod $T\overline N$).  
It is a contradiction (see Theorem \ref{tl}). 
We leave the other cases for the reader's exercise. 
So, there are no non-Gorenstein non-$\mathbb Q$-factorial 
affine toric threefolds with terminal singularities. 
\end{proof}
Therefore, we completed the proof of Theorem \ref{main}. 
\end{proof}

Theorem \ref{main} has a beautiful corollary. 

\begin{cor}[Three-dimensional terminal toric flop]\label{flo}
Let 
$$
\begin{matrix}
X & \dashrightarrow & \ X^+ \\
{\ \ \ \ \ \searrow} & \ &  {\swarrow}\ \ \ \ \\
 \ & W &  
\end{matrix}
$$ 
be a three-dimensional toric flopping diagram such that 
$W$ is affine. 
Assume that $X$ has only terminal singularities. 
Then it is the {\em{simplest flop}}, 
where the simplest flop means the flop described 
in \cite[p.49--p.50]{fulton}.  
\end{cor}
\begin{proof}
By the assumption, $W$ is a non-$\mathbb Q$-factorial 
affine toric threefold with terminal singularities. 
Thus, $X\simeq \Spec \, 
\mathbb C[x,y,z,w]/(xy-zw)$ by Theorem \ref{main}. 
So, the above diagram must be the simplest flop. 
\end{proof}

\section{Three-dimensional $\mathbb Q$-factorial 
terminal toric flips}\label{sec3}
We classify three-dimensional flipping contractions 
from $\mathbb Q$-factorial terminal toric threefolds. The next theorem 
was stated in \cite{kmm} without proof at the end of 
Example 5-2-5. 

\begin{thm}[Three-dimensional $\mathbb Q$-factorial terminal 
toric flips]\label{main1} 
Let $\varphi_R:X(\Delta)\to Y(\Sigma)$ be the contraction morphism 
of an extremal ray $R$ with $K_X\cdot R<0$ of flipping type from 
a toric 
threefold with only $\mathbb Q$-factorial terminal singularities. 
Assume that $Y$ is affine. 
Then we have the following description of the flipping contraction{\em{:}} 

There exist two three-dimensional cones 
\begin{align*}
\tau_4&=\langle v_1, v_2, v_3\rangle\in \Delta, \\
\tau_3&=\langle v_1, v_2, v_4\rangle \in \Delta, 
\end{align*} 
sharing the two-dimensional wall 
$$
w=\langle v_1, v_2\rangle 
$$
such that $[V(w)]\in R$ and that for some $\mathbb Z$-coordinate of $N\simeq 
\mathbb Z^3$, 
\begin{align*}
v_1  &= (1,0,0), & v_2&=(0,1,0), & v_3&=(0,0,1),\\
v_4  &= (a,r-a,-r), 
\end{align*} 
or 
\begin{align*}
v_1  &= (1,0,0), & v_2&=(0,1,0), & v_3&=(0,0,1),\\
v_4  &= (a,1,-r), 
\end{align*} 
where $0<a<r$ and $\gcd (r,a)=1$. 
Therefore, 
$$
\Delta=\{\tau_3, \tau_4, {\text{and their faces}}\}, 
$$ 
and 
$$
\Sigma=\{\langle v_1, v_2, v_3, v_4\rangle, {\text{and its faces}}\}. 
$$
\end{thm}
\begin{proof}
By \cite[Example-Claim 14-2-5]{ma}, 
it is sufficient to prove that the (unique) rational curve that is contracted 
passes through only one singular point of $X$. 
Without loss of generality, we may assume that 
$v_1=(1,0,0)$ and $v_2=(0,1,0)$ since $\langle v_1, v_2\rangle$ is a 
two-dimensional non-singular cone. 
Seeking a contradiction, 
we assume that both $\langle v_1, v_2, v_3\rangle$ and 
$\langle v_1, v_2, v_4\rangle$ are singular. 
By the terminal lemma (\cite[\S 1.6]{oda}), we may assume that 
$v_3=(1,p,q)$, where 
$0<p<q$ and $\gcd (p,q)=1$. 
We note that $q\geq 2$. 
We can write $v_4=av_1+bv_2+c(k, l, -1)$ with 
$0< a<c$, $0< b<c$, $\gcd(a,c)=1$, $\gcd(b,c)=1$, and 
$k, l\in \mathbb Z$. 
In particular, $c\geq 2$. 
We note that we assumed that 
$\langle v_1, v_2, v_4\rangle$ is singular and terminal. 
By the terminal lemma again (see 
\cite[p.36 White's Theorem]{oda}), at least one of $a-1$, $b-1$ and 
$a+b$ is divisible by $c$. 
Therefore, $a=1$, $b=1$, or $a+b=c$. 
We note that $v_1, v_2, v_3$ are on the plane 
$$x+y-\frac{p}{q}z=1.$$ 

\begin{ca}[$a=1$]
In this case, $v_4=(1+ck, b+cl, -c)$. 
We have 
$$\frac{c}{q}v_3+v_4=(1+ck+\frac{c}{q}, b+cl+\frac{p}{q}c, 0). 
$$ 
Thus, we obtain the following three inequalities: 
\begin{equation}\label{1}
1+ck+\frac{c}{q}>0, 
\end{equation}
\begin{equation}\label{2}
b+cl+\frac{p}{q}c>0, 
\end{equation}
and 
\begin{equation}\label{3}
1+ck+b+cl+\frac{p}{q}c<1. 
\end{equation} 
The inequalities (\ref{1}) and (\ref{2}) follow from the condition that $\varphi_R$ 
is small. 
The condition $K_X\cdot R<0$ implies the inequality (\ref{3}). 
By (\ref{2}) and (\ref{3}), we have $k\leq -1$. 
Thus 
$$
0<1+ck+\frac{c}{q}\leq 1-c+\frac{c}{q}\leq 1-\frac{1}{2}c\leq 0
$$ 
by (\ref{1}). It is a contradiction. 
\end{ca}

\begin{ca}[$b=1$]
In this case, $v_4=(a+ck, 1+cl, -c)$. 
We have 
$$\frac{c}{q}v_3+v_4=(a+ck+\frac{c}{q}, 1+cl+\frac{p}{q}c, 0). 
$$ 
Thus, we obtain the following three inequalities: 
\begin{equation}\label{4}
a+ck+\frac{c}{q}>0, 
\end{equation}
\begin{equation}\label{5}
1+cl+\frac{p}{q}c>0, 
\end{equation}
and 
\begin{equation}\label{6}
a+ck+1+cl+\frac{p}{q}c<1. 
\end{equation}
By (\ref{5}) and (\ref{6}), $k\leq -1$. 
So, $k=-1$ by (\ref{4}). 
By (\ref{5}), we know that $l\geq -1$. 
Therefore, $l=0$ or $-1$ by (\ref{6}). 

First, we assume that $l=0$. Then we get 
$$
a-c+\frac{p}{q}c<0
$$
by (\ref{6}) and 
$$
a-c+\frac{c}{q}>0
$$ 
by (\ref{4}). It is a contradiction. 

Next, we assume that $l=-1$. 
Then we obtain 
$$
a-c+\frac{c}{q}>0
$$ 
by (\ref{4}) and 
$$
1-c+\frac{p}{q}c>0
$$ 
by (\ref{5}). 
These two inequalities imply that 
$$
1+a-2c+\frac{p+1}{q}c>0. 
$$
It is a contradiction. 
\end{ca}

\begin{ca}[$a+b=c$]
In this case, $v_4=(a+ck, c-a+cl, -c)$. 
We have 
$$\frac{c}{q}v_3+v_4=(a+ck+\frac{c}{q}, c-a+cl+\frac{p}{q}c, 0). 
$$ 
Thus, we obtain the following three inequalities: 
\begin{equation}\label{7}
a+ck+\frac{c}{q}>0, 
\end{equation}
\begin{equation}\label{8}
c-a+cl+\frac{p}{q}c>0, 
\end{equation}
and 
\begin{equation}\label{9}
a+ck+c-a+cl+\frac{p}{q}c<1. 
\end{equation}
By (\ref{8}) and (\ref{9}), $k\leq -1$. 
So, $k=-1$ by (\ref{7}). 
By (\ref{8}), we have $l\geq -1$. 
Therefore, $l=0$ or $-1$ by (\ref{9}). 

First, we assume that  $l=0$. Then we have 
$$
\frac{p}{q}c<1
$$
by (\ref{9}) and 
$$
a-c+\frac{c}{q}>0
$$ 
by (\ref{7}). 
Thus, 
$$
1>\frac{p}{q}c\geq \frac{c}{q}>c-a\geq 1. 
$$
It is a contradiction. 

Next, we assume that $l=-1$. 
Then we obtain 
$$
a-c+\frac{c}{q}>0
$$ 
by (\ref{7}) and 
$$
-a+\frac{p}{q}c>0
$$ 
by (\ref{8}). 
By adding these two inequalities, we have 
$$
-c+\frac{p+1}{q}c>0. 
$$
It is a contradiction. 
\end{ca}

Therefore, at least one of 
$\langle v_1, v_2, v_3\rangle$ and  
$\langle v_1, v_2, v_4\rangle$ must be non-singular. 
Thus, we have the desired description of $\varphi_R:X\to Y$ 
by \cite[Example-Claim 14-2-5]{ma}. 
\end{proof}
\begin{rem}
The example in \cite[Remark 14-2-7 (ii)]{ma} 
is not true. The cone $\langle v_1, v_2, v_3\rangle$ is not 
terminal. The cone $\langle v_1, v_2, v_3\rangle$ has canonical singularities. 
\end{rem}

\begin{rem}\label{222}
The source space $X$ in Theorem \ref{main1} is always singular. 
\end{rem}
\begin{rem}
In \cite[Example-Claim 14-2-5]{ma}, $X$ is assumed to be 
{\em{complete}}. It is because contraction morphisms of extremal rays are 
constructed only for {\em{complete}} varieties in \cite{reid} and 
\cite[Chapter 14]{ma}. For the details of non-complete toric 
varieties, see 
\cite{fs1}, \cite{fujino}, and \cite{sato}. 
\end{rem}

\section{Main Theorem}\label{sec4}
The following theorem is the main theorem of this 
paper. 

\begin{thm}[cf.~\cite{takano}]\label{main-taka} 
Let $\varphi:X\to Y$ be a small 
proper toric morphism such that $X$ is a three-dimensional 
toric variety with only terminal singularities. 
Let $C\simeq \mathbb P^1$ be an exceptional 
curve of $\varphi$. 
Assume that $-K_X\cdot C>0$. 
Then $C$ does not pass through 
ordinary double points.  
\end{thm}

\begin{proof}
First, we assume that $C$ passes through 
two ordinary double points. 
By taking a small projective 
resolution of $X$, we can assume 
that 
$C$ does not pass through any 
singular points. 
It is a contradiction by Theorem \ref{main1} (see Remark \ref{222}). 

Next, we assume that $C$ passes through only one 
ordinary double points. 
By Theorem \ref{main1}, 
we have the following local description of 
$X$ and $C$: 

There exist lattice points of $N=\mathbb Z^3$  
\begin{align*}
v_1  &= (1,0,0), & v_2&=(0,1,0), & v_3&=(0,0,1),\\
v_5  &= (-1,1,1), & v_6&=(1,-1,1). && 
\end{align*} 
We put  
$$
\Delta_1=\{\langle v_1, v_2, v_3, v_5\rangle, \langle v_1, v_2, 
v_4\rangle, \text{and their faces}\},  
$$ 
and 
$$
\Delta_2=\{\langle v_1, v_2, v_3, v_6\rangle, \langle v_1, v_2, 
v_4\rangle, \text{and their faces}\},  
$$ 
where 
$v_4=(a, r-a, -r)$ or $(a, 1,-r)$ with 
$0<a<r$ and $\gcd(a,r)=1$. 
Then $X=X(\Delta)$, 
where $\Delta =\Delta_1$ 
or $\Delta_2$, and $C$ is $V(\langle v_1, v_2\rangle)\simeq \mathbb P^1$. 
\setcounter{ca}{0}
\begin{ca}
When $v_4=(a, r-a, -r)$ and 
$\Delta=\Delta_1$, 
we have 
$$
v_2=\frac{r}{2r-a}v_5+\frac{r-a}{2r-a}v_1+\frac{1}{2r-a} v_4. 
$$
Therefore, $v_2$ is contained in the cone $\langle v_5, v_1, v_4\rangle$. 
Thus, we can not remove the wall $\langle v_1, v_2\rangle$ from $\Delta$. 
\end{ca}

\begin{ca}
When $v_4=(a, 1, -r)$ and 
$
\Delta=\Delta_1
$, 
we have 
$$
v_2=\frac{r}{r+1}v_5+\frac{r-a}{r+1}v_1+\frac{1}{r+1} v_4. 
$$
Therefore, $v_2$ is contained in the cone $\langle v_5, v_1, v_4\rangle$. 
Thus, we can not remove the wall $\langle v_1, v_2\rangle$ from 
$\Delta$. 
\end{ca}

\begin{ca}
When $v_4=(a, r-a, -r)$ and 
$
\Delta=\Delta_2
$, 
we have 
$$
v_1=\frac{r}{r+a}v_6+\frac{a}{r+a}v_2+\frac{1}{r+a} v_4. 
$$
Therefore, $v_1$ is contained in the cone $\langle v_6, v_2, v_4\rangle$. 
Thus, we can not remove the wall $\langle v_1, v_2\rangle$ from 
$\Delta$. 
\end{ca}

\begin{ca}
When $v_4=(a, 1, -r)$ and 
$\Delta=\Delta_2
$, 
we have 
$$
v_1=\frac{r}{r+a}v_6+\frac{r-1}{r+a}v_2+\frac{1}{r+a} v_4. 
$$
Therefore, $v_1$ is contained in the cone $\langle v_6, v_2, v_4\rangle$. 
Thus, we can not remove the wall $\langle v_1, v_2\rangle$ from 
$\Delta$. 
\end{ca}
Thus, $C$ does not pass through any ordinary double points. 
\end{proof}

\ifx\undefined\bysame
\newcommand{\bysame|{leavemode\hbox to3em{\hrulefill}\,}
\fi


\begin{thebibliography}{KMM}

\bibitem[F1]{fujino}
O.~Fujino, 
Equivariant completions of toric contraction morphisms, 
Tohoku Math. J. {\textbf{58}} (2006), 303--321. 

\bibitem[F2]{spe} 
O.~Fujino, 
Special termination and reduction to pl flips, 63--75 in 
{\em{Flips for 3-folds and 4-folds}}, Oxford University Press (2007).

\bibitem[FS]{fs1} 
O.~Fujino and H.~Sato, 
Introduction to the toric Mori theory, 
Michigan Math. J. {\textbf{52}} (2004), no. 3, 649--665. 

\bibitem[Fl]{fulton}
W.~Fulton, 
{\em{Introduction to toric varieties}}, Annals 
of Mathematics Studies, {\textbf{131}}, The William H. Roever
Lectures in Geometry, Princeton University Press, Princeton, 
NJ, 1993. 

\bibitem[I]{ishida} 
M.~Ishida, 
On the terminal toric singularities of dimension three, 
in {\em{Commutative Algebra}}, Karuizawa, Japan, 1982 
(S.~Goto, ed.), 54--70. 

\bibitem[II]{ishi-iwa} 
M.~Ishida, and N.~Iwashita, Canonical cyclic 
quotient singularities of dimension three, 
{\em{Complex analytic singularities}}, 135--151, 
Adv.~Stud.~Pure Math., {\textbf{8}}, 
North-Holland, Amsterdam, 
1987. 

\bibitem[KMM]{kmm}
Y.~Kawamata, K.~Matsuda, and K.~Matsuki, 
Introduction to the minimal model problem, 
Algebraic geometry, Sendai, 1985, 
283--360, Adv. Stud. Pure Math., {\textbf{10}}, 
North-Holland, Amsterdam, 1987.  

\bibitem[M]{ma} 
K.~Matsuki, {\em{Introduction to the Mori program}}, 
Universitext, Springer-Verlag, New York, 2002. 

\bibitem[O]{oda}
T.~Oda, 
{\em{Convex bodies and algebraic 
geometry. An introduction to the theory of toric varieties}}, 
Translated from the Japanese. Ergebnisse der 
Mathematik und ihrer Grenzgebiete (3) 
[Results in  Mathematics and Related Areas (3)], 
{\textbf{15}}. Springer-Verlag, Berlin, 1988. viii+212 pp.

\bibitem[R]{reid}
M.~Reid, 
Decomposition of toric morphisms, 
Arithmetic and geometry, 
Vol. II, 395--418, Progr. Math. {\textbf{36}}, Birkh\"auser 
Boston, Boston, MA, 1983. 

\bibitem[YPG]{ypg} 
M.~Reid, Young person's guide to canonical 
singularities, {\em{Algebraic geometry, Bowdoin, 1985 
$($Brunswick, Maine, 1985$)$}}, 345--414,  
 Proc. Sympos. Pure Math., {\textbf{46}}, Part 1,  
 Amer. Math. Soc., Providence, RI,  1987. 
 
\bibitem[S]{sato} 
H.~Sato, Combinatorial descriptions of toric extremal 
contractions, Nagoya Math. J. {\textbf{180}} (2005), 111--120. 
 
\bibitem[T]{takano} 
Y.~Takano, 
On flipping contractions of three-dimensional toric varieties with 
non-$\mathbb Q$-factorial terminal singularities (Japanese), 
Master's thesis, Tokyo Metropolitan University (2008).  
\end{thebibliography}
\end{document}